\documentclass[leqno]{amsart}

\usepackage{amsmath,amsfonts,amsthm,amssymb, mathrsfs, mathabx, accents, mathdots, vntex, tikz}
\usepackage[utf8]{inputenc}
\usepackage[french,english]{babel}
\usepackage[bookmarks]{hyperref}
\usepackage{graphics}
\usepackage[all]{xy}
\usepackage{float}

\hypersetup{backref, colorlinks=true}

 \usepackage{graphicx}
\usepackage[all]{xy}
\newtheorem{theorem}{Theorem}[section]

\newtheorem{proposition}[theorem]{Proposition}
\newtheorem{corollary}[theorem]{Corollary}
\newtheorem{remark}[theorem]{{ Remark}}

\newtheorem{example}[theorem]{{ Example}}
\newtheorem{definition}[theorem]{{Definition}}
\newtheorem{definitions}[theorem]{{Definitions}}

\newenvironment{equationth}{\stepcounter{theorem}\begin{equation}}{\end{equation}}
\newenvironment{preuve}{{\em{\noindent \text{{\rm Proof}.} }}}
{\hfill $\Box$} 

\def\C{ \mathbb{C}}
\def\Q{ \mathbb{Q}}
\def\R{ \mathbb{R}}

\def\K{ \mathbb{K}}

\def\dist{ {\rm dist}}

\def\epsilon{\varepsilon} 

\def\rond{\mathaccent"7017}

\begin{document}

\large


\title[]{\small 
Local bi-Lipschitz classification of  
semialgebraic surfaces}
\makeatother

\author{Jean-Paul Brasselet,  Maria Aparecida Soares  Ruas and 
Thuy Nguyen}

\address[Jean-Paul Brasselet]{Institut de Math\'ematiques de Marseille, UMR 7373 du CNRS, 
Aix-Marseille University - Case 907 - Campus de Luminy - 13288 MARSEILLE Cedex 9.}
\email{jean-paul.brasselet@univ-amu.fr}

\address[Thuy Nguyen Thi Bich]{S\~ao Paulo State University, UNESP, S\~ao Jos\'e do Rio Preto, Brazil.}
\email{bich.thuy@unesp.br}

\address[Maria Aparecida Soares Ruas]{Universidade de S\~ao Paulo, Instituto de Ci\^encias Matem\'aticas e de Computa{\c c}\~ao - USP, Avenida Trabalhador S\~ao-Carlense, 400 - Centro, S\~ao Carlos, Brazil.}
\email{maasruas@icmc.usp.br}

\maketitle \thispagestyle{empty}

\begin{abstract} 

We provide 
bi-Lipschitz invariants for finitely determined map germs $f: (\K^n,0) \to (\K^p, 0)$, where $\K = \R$ or $ \C$. The aim of the paper is to provide  
partial answers to the following questions:

Does the bi-Lipschitz type of a map germ $f: (\R^n, 0) \to (\R^p, 0)$ determine the bi-Lipschitz type of 
the link of $f$ and of the double point set of $f$? 
Reciprocally, 
given a map germ $f: (\R^n, 0) \to (\R^p, 0)$, do the bi-Lipschitz types of 
the link of $f$ and of the double point set of $f$ determine the bi-Lipschitz type of the germ 
$f: (\R^n, 0) \to (\R^p, 0)$?

We provide a positive answer to the first question in the case of a finitely determined map germ
$f: (\R^n, 0) \to (\R^p, 0)$ where $p \geq 2n-1$ 
(Theorem \ref{necessarycondition}).

With regard to the second question, for a finitely determined  map germ $f : (\R^2,0) \to (\R^3,0),$ we show (Proposition \ref{demai}) that a complete set of invariants for the bi-Lipschitz
classification of $X_f$ with respect to the inner metric is given by 
the link of $f$, the image of the double point set of $f$ and the polar curve of a generic 
projection into the plane. In particular,  in the homogeneous parametrization case $f: (\R^2, 0) \to (\R^3, 0)$
of corank 1, we do not need the hypothesis on the equivalence of the image of the double point set (Theorem \ref{theo_quasehomogeneo}).
\end{abstract}

\section{INTRODUCTION}

 Two map germs $f, g: (\K^n, 0) \to (\K^p, 0)$, where $\K = \C$ or $ \R$, are $\mathcal{A}$-equivalent
 if there exist diffeomorphism map germs $h: (\K^n, 0) \to (\K^n, 0)$ and $k: (\K^p, 0) \to (\K^p, 0)$ such that  the following diagram is commutative

\def\banane{\ar@{}[rd]|{\circlearrowleft }} 
\begin{equation}\label{diag1}
\xymatrix @C=2cm{ 
(\K^n, 0)    \ar[r]^{f}       \ar[d]_{h}       \banane      &  (\K^p, 0)       \ar[d]^{k} \\ 
(\K^n, 0)     \ar[r]_{g}     &  (\K^p, 0),
}
\end{equation}
that means $g = k \circ f \circ h^{-1}$. 
When $h$ and $k$ are homeomorphisms, then $f$ and $g$ are called topologically $\mathcal{A}$-equivalent, or ${\mathcal{C}^0}$-$\mathcal{A}$-equivalent. 
Two map germs  $f$ and  $g$ are {\it bi-Lipschitz  ${\mathcal A}$-equivalent} if   $h$ and $k$ in the above diagram are bi-Lipschitz homeomorphisms. 

Note that bi-Lipschitz equivalence implies topological equivalence. 

 A map germ $f: (\K^n, 0) \to (\K^p, 0)$ is ${\mathcal A}$-finitely determined  if there exists a positive integer $r$ such that for any $g$ with equality of the $r$-jets $j^r g (0) = j^r f (0)$, 
 then $f$ and $g$ are $\mathcal{A}$-equivalent.
 In this case, we say that $f$ is $r$-${\mathcal A}$-finitely determined.

In this paper, we will consider only  ${\mathcal A}$-equivalence classes and we will omit the 
symbol ${\mathcal A}$. 
Moreover, we consider the case $p\geq 2n-1$ and $n \geq 2$.

\medskip

In these dimensions, 
the natural definition of double point set of $f$, that we will consider in this paper, denoted by 
${D}(f)$,  is geometrically, defined by 
$$ D(f) = \{ x\in U : \exists  x' \ne x, f(x) = f(x') \} \cup \Sigma(f)$$
where $U$ is a sufficiently small neighbourhood of $0$ and  $\Sigma(f)$ is the singular set of $f$.

We also consider the lifting of $D(f)$ in $U\times U$, denoted by $\widetilde{D}(f)$ given by the pairs $(x,x')$ such that either $f(x) = f(x')$ with $x\ne x'$ or $x=x'$ with $x\in \Sigma(f)$. 


An alternative definition of $D(f)$ has been provided by D. Mond \cite{Mond1} in terms
of ideals. It has been used by various authors, for example \cite{MNB,Thon1,Ot}. 

 \medskip
 
Consider a representative of a finitely determined map germ $f: U \to (\K^p, 0)$, where $U \subset \mathbb K^n$ and $n\leq p.$ Let $X_f$ denote the image  $f(U)$ of $f$. The intersection of $X_f$ with a sufficiently  small sphere $S_{\epsilon}^m$  centered at the origin in $\K^p$  is called  the {\it link} of the map germ $f$, and denoted by $link(f)$. Here, $m=p-1$ in the real case and $m=2p-1$ in
the complex situation. 

According to Mather-Gaffney geometric criteria, \cite{Wall}, when  $p \geq 2n-1$ and $n\geq 2$, 
the singularities of finitely determined map germ 
$f$ are isolated.  
 When $p$ is greater than $2n$, in a punctored neighborhood of the origin, 
$f$ is an embedding  \cite{Gaf,Wall}. 
Note that, according to Mather  \cite{Mat}, the dimensions 
$(n,p)$ satisfying $p \geq 2n-1$ are nice dimensions in the sense of Mather 
(see for example \cite[Table 1, page 181]{DPW}). 

In this work, we consider the following problem: 
Let $f, g: (\K^n, 0) \to (\K^p, 0)$, where $p\geq 2n-1$,
    be two finitely determined map germs. If $f$ and $g$ are bi-Lipschitz equivalent, are 
the links $link(f)$ and $link(g)$ and 
the images of the double point sets $f(D(f))$ and $g(D(g))$  bi-Lipschitz equivalent sets ?
Does the converse hold ? If not, what is the complete set of invariance
for the  bi-Lipschitz equivalence of $f$ and $g$?


\smallskip 

Using equisingularity in family and in the complex case $\K=\C$, 
a counter-example to the converse has been provided by O. N. Silva \cite{Ot, Ot1}. In this example 
the real links are topologically equivalent  and the images $f(D(f))$ of the double point sets 
are bi-Lipschitz equivalent, but the germs are not bi-Lipschitz equivalent. 

In the complex case, 
 it would be interesting to characterize classes of maps and conditions 
 for which the converse is true.

In this paper, we provide some elements for solving the problem and its converse.

\medskip 

The paper has three parts:

\medskip

I)  The first part presented in the section \ref{Bi-Lipschitz_Invariants} provides some bi-Lipschitz invariants, for both cases $\K = \R$ and $\K = \C$,  of the double point sets and of the real links  for the general case $n \leq p$ and for the case  $p \geq 2n-1$ (Proposition  \ref{invdoublepoints}). 
In particular, in the real case,  we show 
 that if $f$ and $g$ are bi-Lipschitz equivalent, then  
the links $link(f)$ and $link(g)$ and 
the images of the double point sets $f(D(f))$ and $g(D(g))$ are bi-Lipschitz equivalent (Theorem \ref{necessarycondition}).

\medskip 

 II)  In the second part, presented in the section  \ref{Bi-Lipschitz_classification}, 
we show that, in general, for a 
finitely determined map germ $f : (\R^2,0) \to (\R^3,0),$ a complete set of invariants for the bi-Lipschitz
classification of $X_f$ with respect to the inner metric (see Definition \ref{def4.1}) is given by 
the link of $f$, the image of the double point set of $f$ and the polar curve of a generic 
projection into the plane. 

For this weaker equivalence relation, 
the proposed problem is completely solved (Proposition \ref{demai}). The main tools used in this section are the 
Birbrair construction of H\"older Complexes (Definitions \ref{contact_Holder} and \ref{geometrie}) and  Birbrair Classification Theorem \ref{Clas_theo} (\cite{Lev}). 

\medskip 

III) The third part is presented in the section  \ref{BiLipschitzforquasihomogenousSurfaces}. 
We prove that for a corank 1 finitely determined map germ with a homogeneous parametrization $f: (\R^2, 0) \to (\R^3, 0),$ 
we do not need the hypothesis on the equivalence of the image of the double point locus, {\it i.e} in this case 
 the bi-Lipschitz type of 
the link of $f$ and the polar curve of a generic 
projection into the plane determine the bi-Lipschitz type of $f$ (Theorem \ref{theo_quasehomogeneo}). 
This result provides a class of quasi-homogeneous surfaces for which the bi-Lipschitz type of the 
link and the polar curve determine the bi-Lipschitz type of the surface. 
Finally, we  apply our results to relate the $C^{0}- \mathcal A$ classes of finitely determined map germs $f$ of corank 1 with  homogeneous  parametrization 
 and the inner bi-Lipschitz type of $X_f$ (Proposition \ref{Proposition 5.4}).

\medskip 

Results on the related problem of $C^{0}$-$\mathcal A$-equivalence of finitely determined germs were given by J. Nu\~{n}o Ballesteros and W. Marar in \cite{Thon1} (case $n=2, p=3$) and by R. Mendes and J. Nu\~{n}o Ballesteros in \cite{NR} (case $n=2$ and $p=4$.) 
We also refer to recent results by W. Neuman and A. Pichon \cite{NP}, 
Paunescu  and Parusi\'nski \cite{PP1} \cite{PP2} and R. Mendes and J.E. Sampaio \cite{MaSa}. 

\section{Some well-known results}
We present in this section some well-known results that will be used in the next sections of the paper.

\begin{theorem} [Mather-Gaffney Geometric criterion, \cite{Wall}] \label{criteriongeometric}
{\rm A map germ $f: (\C^n, 0) \to (\C^p, 0)$ is finitely determined if and only if, for every representative $f$ of the map germ, there exist a neighbourhood $U$ of 0 in $\C^n$ and a  neighbourhood $V$ of 0  in $\C^p$, with $f(U) \subset V$, such that for all $y \in V \setminus \{0\}$, the set $S = f^{-1}(y) \cap \Sigma(f)$ is finite and $f: (\C^n, S) \to (\C^p, y)$ is stable, where $\Sigma(f)$ is the set of critical points of $f$. 

}
\end{theorem}
 
 Let $f: (\C^2, 0) \to (\C^3, 0)$ 
be a finitely determined map germ. From the classical result of Whitney \cite{Whitney}, 
we know that the stable singularities in these dimensions are transverse double points, triple points and cross-caps. In this case,  Theorem \ref{criteriongeometric} says that 
$f: (\C^2, 0) \to (\C^3, 0)$ is ${\mathcal A}$-finitely determined if and only if for every representative $f$, 
there exists a neighbourhood $U$ of 0, such that the only singularities of 
$f(U) \setminus \{ 0 \}$ are transverse double points. 

Notice that,  in the real case, {\it i.e.} in the case  of map germs $f: (\R^n, 0) \to (\R^p, 0)$, the converse implication  of Theorem \ref{criteriongeometric} does not hold.
 However, we have

\begin{proposition} \label{criteriongeometric_real}
{\rm If a map germ $f: (\R^2, 0) \to (\R^3, 0)$ is finitely determined then there exist a neighbourhood $U$ of 0 in $\R^2$ and a  neighbourhood $V$ of 0  in $\R^3$, with $f(U) \subset V$, such that the only singularities of 
$f(U) \setminus \{ 0 \}$ are transverse double points.

}
\end{proposition}

In these dimensions, the $C^{0}-\mathcal A$ classification of finitely determined germs has been considered by 
 Nu\~no Ballesteros and Marar in \cite{Thon1}. 

The link of a finitely determined map germ is a closed curve which is a  {\it stable immersion}  $\gamma : S^1 \to S^2.$  This means 
that the singular points of its image $f(S^1)$ are  transverse double points (also called {\it crossings}). In \cite{Thon1}, the authors call the link of $f$ as the associated
{\it doodle } of $f.$ The following result was given in \cite{Thon1}.

\begin{theorem}{\rm (\cite{Thon1}, Corollary 3.4) \label{theorem-thon-juan-2} Two finitely determined map germs $ f, g: (\mathbb R^2,0) \to (\mathbb R^3,0)$ are topologically equivalent if and only if their associated doodles are topologically equivalent.}
	\end{theorem}

To describe the topology of  the link,  one can associate a {\it Gauss word} to the doodle (see \cite{Thon1} for more details).
	
	\begin{definition}{\rm  (\cite{Thon1} Definition 3.6) 
	Let $\gamma : S^1 \to S^2$ be a doodle with $r$ crossings. Choose $r$ letters
			$\{a_1, a_2, \ldots, a_r\}$ to label them. Fix orientations on $S^1$ and $S^2$, and choose a base point $z_0$ in $S^1.$
			Now, consider a permutation 
			$$\sigma:  \{1, 2,\ldots, 2r\} \to \{a_1, a_2, \ldots, a_r, a_1^{-1},\ldots, a_r^{-1}\},$$ constructed as follows:
			Let $z_1, \ldots, z_{2r}$ denote the source double points ordered as $z_0 \leq z_1 < z_2< \ldots < z_{2r}.$ Assume that
			$\gamma(z_i)=\gamma(z_j)=a_k,\, i < j.$ Then, write $\sigma(i)=a_k,\, \sigma(j)=a_k^{-1}$ if the pair of vectors $(\gamma'(z_i), \gamma'(z_j))$ are positively oriented 	in $S^2, $ or $\sigma(i)=a_k^{-1}$ and $\sigma(j)=a_k, $ otherwise.} 
			The sequence	$\sigma(1) \ldots \sigma(2r)$ is called the (signed) Gauss word of the 
doodle $\gamma$.
		\end{definition}

	The Gauss word is not uniquely determined, 
	since it depends on the labels $\{a_1, a_2, \ldots, a_r\}$, the chosen orientations on both $S^1$ and 
	$S^2$ and on the base point $z_0 \in S^1$. In \cite{Thon1}  the authors explicit equivalences
	between Gauss words in such a way that up to equivalences, the Gauss word is well defined and the following holds
	
\begin{corollary}{\rm (\cite{Thon1}, Corollary 3.8)} \label{theorem-thon-juan} 
Two finitely determined map germs $ f, g: (\mathbb R^2, 0) \to (\mathbb R^3,0)$ are topologically equivalent if and only if the Gauss words of their associate doodles are  equivalent.
	\end{corollary}

The following well-known theorem of Fukuda \cite{Fukuda} will be also used in the next sections of the paper:
\begin{theorem} \cite{Fukuda} \label{theoremFukuda}
{\rm Suppose that $n \leq p$. Then given a semialgebraic subset $W$ of $J^r(n, p)$, 
there exist an integer $s \geq r$, depending only on $n, p$ and $r$, 
and a closed semialgebraic subset $\Sigma_W$ of $(\pi_r^s)^{-1}(W)$, where $\pi_r^s: J^s(n, p) \to J^r(n, p)$ is the canonical projection, 
having codimension greater than 1 such that for any $C^{\infty}$-mapping $f: \R^n \to \R^p$ 
with $j^sf(0)$ belonging to $(\pi_r^s)^{-1}(W) \setminus \Sigma_W$, there exists a positive number $\epsilon_0$ such that for any number $\epsilon$ with $0 < \epsilon \leq \epsilon_0$ we have
\begin{enumerate}
\item $\tilde{S}_\epsilon^{n-1} = f^{-1}(S_\epsilon^{p-1})$ is a homotopy ($n-1$)-sphere which, 
if $n \neq 4, 5$, is diffeomorphic to the natural ($n-1$)-sphere $S^{n-1}$,

\item the restricted mapping $f\vert_{\tilde{S}_\epsilon^{n-1}}: \tilde{S}_\epsilon^{n-1} \to {S}_\epsilon^{p-1}$ is to\-po\-lo\-gi\-cally stable. Moreover,  $f\vert_{\tilde{S}_\epsilon^{n-1}}$ is $C^{\infty}$-stable if ($n, p$) are nice dimensions,

\item denoting by $\tilde{B}_\epsilon^{n} = f^{-1}({B}_\epsilon^{p})$ the inverse image of the 
$p$-dimensional ball of radius  $\epsilon$ centered at $0$, 
the restricted mapping $f\vert_{\tilde{B}_\epsilon^{n}} : \tilde{B}_\epsilon^{n} \setminus \{ 0\} \to {B}_\epsilon^{p} \setminus \{ 0\}$ is proper, to\-po\-lo\-gi\-cally stable ($C^{\infty}$ stable if $(n,p)$ are nice dimensions) and topologically e\-qui\-va\-lent ($C^{\infty}$ equivalent if $(n,p)$ are nice dimensions) to the product mapping 
$$(f\vert_{\tilde{S}_\epsilon^{n-1}}) \times id_{(0, \epsilon)} : \tilde{S}_\epsilon^{n-1} \times (0, \epsilon) \to S_{\epsilon}^{p-1} \times (0, \epsilon)$$
defined by $(x,t) \mapsto (f(x), t)$ and 

\item consequently, $f\vert_{\tilde{B}_\epsilon^{n}} : \tilde{B}_\epsilon^{n} \to {B}_\epsilon^{p}$ is topologically equivalent to the cone 
$$C(f\vert_{\tilde{S}_\epsilon^{n-1}}): \tilde{S}_\epsilon^{n-1} \times [0, \epsilon) / \tilde{S}_\epsilon^{n-1} \times \{0\} \to {S}_\epsilon^{p-1} \times [0, \epsilon) / {S}_\epsilon^{p-1} \times \{0\}$$
of the stable mapping $f\vert_{\tilde{S}_\epsilon^{n-1}}: \tilde{S}_\epsilon^{n-1} \to {S}_\epsilon^{p-1}$ defined by $C(f\vert_{\tilde{S}_\epsilon^{n-1}})(x,t) = (f(x), t).$
\end{enumerate} 
}
\end{theorem}

\section{Some bi-Lipschitz Invariants} \label{Bi-Lipschitz_Invariants}

In this section, we consider the bi-Lipschitz equivalence in the ambient space. This means that two subsets $X, Y \subset \K^p$, where $\K = \R$ or $\K = \C$, are bi-Lipschitz equivalent if there exists a homeomorphism in the ambient space $h: \K^p \to \K^p$ such that $h(X) = Y$. 

\begin{proposition} \label{invdoublepoints}
 Let  $f, g: (\K^n, 0) \to (\K^p, 0)$ be two finite analytic map germs,   where $\K = \R$ or $\K = \C$ and 
$p \geq 2n-1$. 
If $f$ and $g$ are bi-Lipschitz equivalent, then 

\begin{enumerate}
\item $D(f)$ and $D(g)$ are  bi-Lipschitz equivalent;
\item $\widetilde{D}(f)$ and $\widetilde{D}(g)$ are  bi-Lipschitz equivalent; 
\item $f(D(f))$ and $g(D(g))$ are bi-Lipschitz equivalent.
\end{enumerate}

\end{proposition}

\begin{preuve} 
Assume that $f$ and $g$ are bi-Lipschitz equivalent, then there exist bi-Lipschitz  homeomorphisms 
$$h: (\K^n, 0) \to (\K^n, 0), \quad k: (\K^p, 0) \to (\K^p, 0)$$
such that $f = k^{-1} \circ g \circ h.$

\medskip 

(1) We will prove  that $h(D(f)) = D(g)$. 
 Take $x \in D(f)$, then there are two cases:
 
\medskip 

a) There exists $x' \in \K^n$ such that $x \neq x'$ and $f(x) = f(x')$. 
Since $f = k^{-1} \circ g \circ h$, then $k \circ f =g \circ h$, one has 
$$k(f(x)) = g(h(x)), \quad k(f(x')) = g(h(x')).$$
Since $f(x) = f(x')$, then $k(f(x)) = k(f(x'))$, hence $g(h(x)) = g(h(x'))$. 
Moreover, since $h$ is bijective, then $h(x) \neq h(x')$.  
Consequently, $h(x)$ is a double point of $g$, or $h(x) \in D(g)$. 

\medskip  

b) If (a) does not occur and $x\in D(f)$, then $x\in \Sigma (f)$. In this case, by 
\cite[Lemma 4.6]{NRT} the singular sets 
$\Sigma (f)$ and $\Sigma (g)$ are bi-Lipschitz equivalent. 
The proof for ${\mathcal A}$--equivalence follows easily. 

We conclude that  $h(D(f)) \subset D(g)$.
\medskip 

(2) 
Remember that $\widetilde{D}(f)$ is the lifting of $D(f)$ in $U\times U$.
The fact ``$\widetilde{D}(f)$ and $\widetilde{D}(g)$ are bi-Lipschitz equivalent'' comes directly from (1). In fact, in this case $(h \times h) (\widetilde{D}(f)) = \widetilde{D}(g) $, where $(h \times h)(x, y) = (h(x), h(y))$, with $(x,y) \in \widetilde{D}(f)$.
\medskip 

(3) We prove now that $f(D(f))$ and $g(D(g))$ are bi-Lipschitz equivalent. 
In fact, we prove that $k(f(D(f))) = g(D(g))$. 

Take $y \in f(D(f))$ and $z = k(y)$, we want to show that $z\in g(D(g))$. 

Let $y=f(x)$  there are two possibilities (for $x$): 

(i) either  there exist $x'\ne x$ in $f^{-1} (y)$.
Then $h(x) \ne h(x')$ and  $g(h(x)) = g(h(x')) = z$. Then $z \in g(D(g))$; 

(ii) or $x\in \Sigma (f)$. Then  $h(x) \in \Sigma (g) \subset D(g)$. 
Since $k \circ f =g \circ h$ then $k(f(x)) = g(h(x))$, therefore $k(y) = g(h(x))$. 
Consequently, $k(y) \in g(D(g))$ and  hence $k(f(D(f))) \subset g(D(g))$. 
We proceed similarly to prove that $g(D(g))\subset k(f(D(f))),$ replacing $k$ by $k^{-1}.$
\end{preuve}
\medskip

 If $f$ is finitely determined and $p=2n-1,$ 
  the  image of the set of double points of $f$ 
 is a curve, which can be embedded in the target space $\K^{2n-1};$  if $p>2n-1,$ the singularity of the image $X_f = f(U)$ is isolated 
\cite{Wall, Whitney}.

\begin{proposition}  \label{inv_link} 
 Let $f, g: (\R^n, 0) \to (\R^p, 0)$ be two finitely determined map germs, 
where  $p\geq 2n-1$. If $f$ and $g$ are bi-Lipschitz equivalent, then $link(f)$ and $link(g)$ are  bi-Lipschitz equivalent. 
\end{proposition}

\begin{preuve}
Let $f: (\R^n, 0) \to (\R^p, 0)$, with $p\geq 2n-1$  be a finitely determined map germ, then there exists an open subset $U$ of $0$ in $\R^n$ and an  open subset $W$ of $0$ in $\R^p$ such that 
$$f \vert_{U \setminus \{ 0\}}: U \setminus \{ 0\} \to W$$
is an immersion whose singularities are at most transverse double points. 
 By Theorem \ref{theoremFukuda}, since $p\geq 2n-1$, then for $\epsilon$ enough small, $f^{-1}(S_{\epsilon}^{p-1})$ is homeomorphic to a sphere $S_\epsilon^{n-1}$ for some homeomorphism 
 $\varphi : f^{-1}(S_{\epsilon}^{p-1}) \to S_\epsilon^{n-1}$. Then 
$$f \vert_{S_\epsilon^{n-1}} \circ \varphi^{-1}  : S_\epsilon^{n-1} \to S_\epsilon^{p-1}$$
is topologically stable. 
The same thing happens with $g$. 
 Furthermore, in this case, $f \vert _{U \setminus \{ 0 \}}$ is transverse to $S_\epsilon^{p-1}$. It follows that  the inverse image $f^{-1}(S_\epsilon^{p-1})$ is diffeomorphisc  to the sphere $S_\epsilon^{n-1}$. 
Then, the images of the maps 
$$f\vert_{f^{-1}(S_\epsilon^{p-1})} :f^{-1}(S_\epsilon^{p-1}) \to S_\epsilon^{p-1} \quad \text{ and } \quad 
 g\vert_{g^{-1}(S_\epsilon^{p-1})} :g^{-1}(S_\epsilon^{p-1}) \to S_\epsilon^{p-1}$$
are respectively, $link(f)$ and $link(g)$. 
Now, the bi-Lipschitz equivalence of $link(f)$ and $link(g)$ follows from the bi-Lipschitz equivalence of $f$ and $g$. 
\end{preuve}

\medskip 

From Propositions \ref{invdoublepoints} and \ref{inv_link}, 
we have the following theorem: 
\begin{theorem} \label{necessarycondition} 
{Let $f, g: (\R^n, 0) \to (\R^p, 0)$  be two finitely determined map germs, 
where  $p\geq 2n-1$. 
If  $f$ and $g$ are  bi-Lipschitz equivalent, then $link(f)$ and $link(g)$ are bi-Lipschitz 
equivalent and the images of  the double point
sets
$f(D(f))$ and $g(D(g))$ are bi-Lipschitz equivalent.}
\end{theorem}

To prove the converse, we need to consider another invariant, the Lipschitz type of  the polar curve
$\Gamma^1(f).$ In \cite[Example 5.2]{Ot} Ruas and Silva show that the family of complex polynomial mappings 
$$f(x,y,t) = (x^2+txy, xy^2+x^2y+y^3, x^5+y^5)$$ is topologically trivial but it is not Whitney equisingular, so that it cannot be bi-Lipschitz trivial. Moreover, for $t=0$ and $t=\bar t,\, \bar t \neq 0,$ the following holds:
\begin{itemize}
	\item $link(f_0)\, \simeq_{bi-Lip} \, link{f_{\bar t}}$
	\item $f_0(D(f_0))\, \simeq_{bi-Lip} \, f_{\bar t}(D(f_{\bar t} ))$
	\end{itemize}
However the polar sets $\Gamma^1(f_0)$ and $\Gamma^1(f_{\bar t})$ are not bi-Lipschitz equivalent.

\section{Bi-Lipschitz classification in the case $(n,p) = (2,3)$.} \label{Bi-Lipschitz_classification}

We start this section with the following result, which will  be useful in the next sections:
\begin{proposition}\label{link-cond}
	Let $f,g: (\mathbb R^2,0) \to (\mathbb R^3,0)$ be  finitely determined map-germs. Then the following are equivalent
	\begin{itemize}
		\item[(a)] $link(f)$ and $link(g)$ are bi-Lipschitz equivalent;
			\item[(b)] $link(f)$ and $link(g)$ are topologically equivalent;
			\item [(c)] the Gauss words associated to $link(f)$ and $link(g)$ are equivalent.
	\end{itemize}
\end{proposition}

\begin{preuve}
Let $S^2_{\delta}$ be a sphere with center at the origin in $\mathbb R^3$ and radius $\delta.$	
Taking representatives of $f$ and $g$ defined on some small neighbourhood $U$ of the origin, the finite determinacy of the germs implies that  the set $f^{-1}(f(U)\cap S_{\delta}^2)=\tilde S^1$ is diffeomorphic to $S^1$
	and the map $f|_{\tilde S^1}: {\tilde S^1} \to S_{\delta}^2$ is an immersion with normal crossings, hence it is a stable mapping. The same holds to $g.$ In these dimensions, two stable mappings are $\mathcal A$-equivalent if and only if they are topologically equivalent (see \cite{DPW}). Since the source, $\tilde S^1$ is compact, they are also bi-Lipschitz  $\mathcal A$-equivalent. Then it follows that (a) and (b) are equivalent.
	The equivalence between (a) and (c) follows from Theorem 
	\ref{theorem-thon-juan-2} and 
	Corollary \ref{theorem-thon-juan}.
		\end{preuve}

\subsection{H\"older Complexes}

We recall here some definitions  which will be useful for our results later on. 

There are two natural metrics defined on 
the spaces $X\subset \R^p$:
the {\it outer metric} or {\it euclidean metric} $$d_{\rm out}(x, y) = \Vert x-y \Vert$$
which is the induced Euclidean metric on $X$ and the {\it inner metric} or {\it length metric} 
$$d_{\rm in}(x, y) = \inf_{\lambda \in \Lambda(x,y)}  l(\lambda)$$
where $\Lambda(x,y)$ is the set of rectifiable arcs $\lambda : [0,1] \to X$ with $\lambda(0) = x$ and $\lambda(1) = y$ and $l(\lambda)$ is the length of $\lambda$.
It is clear that the condition $d_{\rm out}(x, y) \le d_{\rm in}(x, y)$ holds but the converse does not hold in general.
\begin{definition}\label{def4.1}
{\rm 
We say that $X$ is {\it Lipschitz normally embedded} (LNE) if there exists $k > 0$ such that for all $x, y \in X$, we have
$$d_{\rm in}(x, y) \le k \, d_{\rm out}(x, y).$$}
\end{definition}

 Notice that as the germ $f:(\R^2,0) \to (\R^p,0)$, with $p\ge 3$,   is finitely determined, then we can consider $f$ as the germ of a polynomial map. 
Let $U$ be a neighbourhood of $0$ in $\R^2$, then $f(U)$ is a semialgebraic surface of $\R^p$.

In this section, we consider the following equivalence relation between two semialgebraic sets:
\begin{definition}
{\rm 
Let $X$ and $Y$ be  two semialgebraic sets in $\R^p$
and $d_X$ and $d_Y$ be chosen me\-trics in $X$ and $Y$, respectively. 
We say that $X$ and $Y$ are {\it abstract bi-Lipschitz equivalent}  if 
there exists a bi-Lipschitz homeomorphism 
$h : (X, d_X ) \to (Y, d_Y )$ such that $h(X) = Y$. When the homeomorphism $h$ is semialgebraic, we say that $X$ and $Y$ are (abstract) semialgebraically bi-Lipschitz equivalent.
}
\end{definition}

When $d_X$ and $d_Y$ are the inner metrics, 
we say that $X$ and $Y$ are (abstract) bi-Lipschitz inner equivalent. 
Similarly, if $d_X$ and $d_Y$ are the outer metrics, 
we say that $X$ and $Y$ are (abstract) bi-Lipschitz outer equivalent. 
Notice that, in these cases, the homeomorphisms $h : (X, d_X ) \to (Y, d_Y )$ are not necessarily 
defined in the whole ambient space. 

\begin{definition}[See \cite{Cidinha}] \label{contact} 
{\rm 
Let $A$ and $B$ be two arcs in $f(U)$, we define the {\it contact order} between $A$ and $B$ as
$$K(A, B) = ord_r \left(\dist(A \cap S(0, r), B \cap S(0, r))\right),$$
 where $S(0,r)$ is the sphere of center $0$ and radius $r$ in $\R^3$.}
\end{definition}

 Lev Birbrair \cite{Lev} introduced a construction called H\"older Complex in order to study bi-Lipschitz classification of semialgebraic surfaces. We provide some definitions and results of this construction. Notice that this study is performed on the surfaces, and the metric is the inner metric.  Birbrair's construction is summarized as follows: 

Let $\Gamma$ be a finite graph. We denote by 
$E_{\Gamma}$  the set of edges and by $V_{\Gamma}$ the set of vertices of $\Gamma$.
 
\begin{definitions} \label{contact_Holder} 
{\rm 
A {\it H\"older Complex} is a pair $(\Gamma, \beta)$, where $\beta: E_\Gamma \to \Q$  such that $\beta(g) \geq 1$, for every $g \in E_\Gamma$.
Two H\"older Complexes $(\Gamma_1, \beta_1)$ and  $(\Gamma_2, \beta_2)$ are called {\it combinatorially equivalent} if there exists a graph isomorphism $i: \Gamma_1 \to \Gamma_2$
 such that, for every  $g \in E_{\Gamma_1}$, we have $\beta_2(i(g)) = \beta_1(g)$.}
\end{definitions}

\begin{definitions} \label{contact_Holder}
{\rm The {\it standard $\beta$-H\"older triangle} $T_\beta$, with $\beta \in [1, \infty[ \, \cap \, \Q$, is the semialgebraic subset of $\R^2$ defined by
$$T_\beta = \{ (x, y) \in \R^2: 0 \leq y \leq x^\beta, 0 \leq x \leq 1\}.$$
A (semialgebraic) subset $X$ of $\R^3$ is called a {\it $\beta$-H\"older triangle} with the principal vertex $a \in X$ if the germ $(X, a)$ is (semialgebraically) outer bi-Lipschitz equivalent to the germ $(T_\beta, 0)$. 
}
\end{definitions}

\begin{figure}[H]
\begin{tikzpicture}   

\draw  (-0.2, -0.05) cos (0,0) sin (1,0.5) 
cos (2,0) 
;

\draw (2,0) sin (2.2, 0.05);
 \draw (2,0) sin (2.2, -0.05);
\draw (0, 0) parabola (0.5,-2);
\draw (0.5,-2) sin (2,0);
\node at (0.8,-0.3) {$(X,a)$};
\node at (0,-1) {$\gamma_1$};
\node at (1.5,-1) {$\gamma_2$};
\node at (0.5,-2)[below] {$a$};

\node at (5,-1) {$\simeq$};

\draw[->] (7,-2) -- (7,0);
\node at (7,0)[above] {$y$};
\node at (10,-2)[right] {$x$};
\draw (7,-2) parabola (8.5,-0.2);
\draw  (8.5,-0.2)--(8.5, -2);

\fill[fill=gray!45] (7,-2) parabola (8.5,-0.2) -- (8.5, -2) --  (7,-2) ;
\draw[->] (7,-2) -- (10,-2);

\node at (8.5, -2)[below] {1};

\node at (7,-2)[below] {0};
\node at (8.2, -1.5) {$T_\beta$};
\node at (8.5,0){$y=x^\beta$};
\end{tikzpicture} 
\caption{$\beta$-H\"older triangle and the standard one.} \label{figure_contact}
\end{figure}
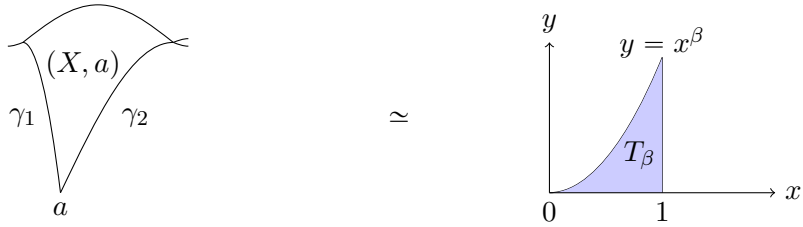

\begin{remark}
{\rm Let $X$ be a $\beta$-H\"older triangle then the contact order between the two branches $\gamma_1$ and $\gamma_2$ is equal to $\beta$ (see Figure \ref{figure_contact}).}
\end{remark}

\begin{theorem} [Theorem 4.1, \cite{LA}] \label{LA} 

{\rm 
Let  $(X,0)$ and $(Y,0)$  be semialgebraic curves with branches $X_1,\ldots ,X_l$ and $Y_1,\ldots ,Y_s$. Then $(X, 0)$ is outer bi-Lipschitz equivalent to $(Y, 0)$ if and only if $l = s$ and there is a permutation $\sigma$ of $\{1,\ldots,l\}$ such that
$$K(X_i,X_j) = K(Y_{\sigma(i)},Y_{\sigma(j)} ) \quad i,j \in \{1,2,\ldots,l \}.$$
}
\end{theorem}

\begin{definition}\label{geometrie}
{\rm Let $(\Gamma, \beta)$ be a H\"older Complex. 
A set $X \subset \R^3$ is called a (semialgebraic) {\it  Geometric H\"older Complex} corresponding to $(\Gamma, \beta)$  with the principal vertex $a \in X$ if:
\begin{enumerate}
\item There exists a homeomorphism $\psi: C\Gamma \to X$, where $C\Gamma$ is the topological cone over $\Gamma$.
\item Let $\alpha \subset C \Gamma$ 
be the vertex of $C\Gamma$,  then $\psi(\alpha) =a$.
\item For each $g \in E_\Gamma$, the set $\psi(Cg)$ is a (semialgebraic) 
$\beta(g)$-H\"older triangle with the principal vertex $a$, 
where $Cg \subset C\Gamma$ is the subcone over $g$. 
\end{enumerate}
}
\end{definition}

\begin{theorem} [Theorem 6.1, \cite{Lev}] \label{theorem4.7}
{\rm 
Let $X \subset \R^p$, with $p\ge 3$,  be a two-dimensional closed semialgebraic set and let $a \in X$. Then there exist a number $\delta >0$ 
and a H\"older Complex $(\Gamma, \beta)$ 
such that $B(a, \delta) \cap X$ is a semialgebraic Geometric H\"older Complex corresponding to $(\Gamma, \beta)$ with the principal vertex $a$, where  $B(a, \delta)$ is the closed ball centered at $a$ of radius $\delta$.   
}
\end{theorem}

\begin{definition}
{\rm We say that $b$ is a {\it non-critical vertex} of $\Gamma$ 
if it is incident with exactly two different edges $g_1$ and $g_2$ and these edges connect two different vertices $b_1$ and $b_2$ with $b$.

If this vertex $b$ is connected by $g_1$ and $g_2$ with only one other vertex $b'$, we say that $b$ is a {\it loop vertex}.

The other vertices of $\Gamma$ (which are neither non-critical nor loop) are called {\it critical vertices} of $\Gamma$.
}
\end{definition}

Given a Geometric H\"older Complex, there is a simplification process described by Lev Birbrair (Theorem 7.3, \cite{Lev}), allowing to assume that every vertex in $V_\Gamma$ is either a critical vertex or a loop vertex.  The resulting Geometric H\"older Complex is called a {\it Canonical H\"older Complex} of $X$ at $a$. Two simplifications of the same H\"older Complex are combinatorially equivalent. 

\begin{figure}[H]
\scalebox{0.50}{\includegraphics{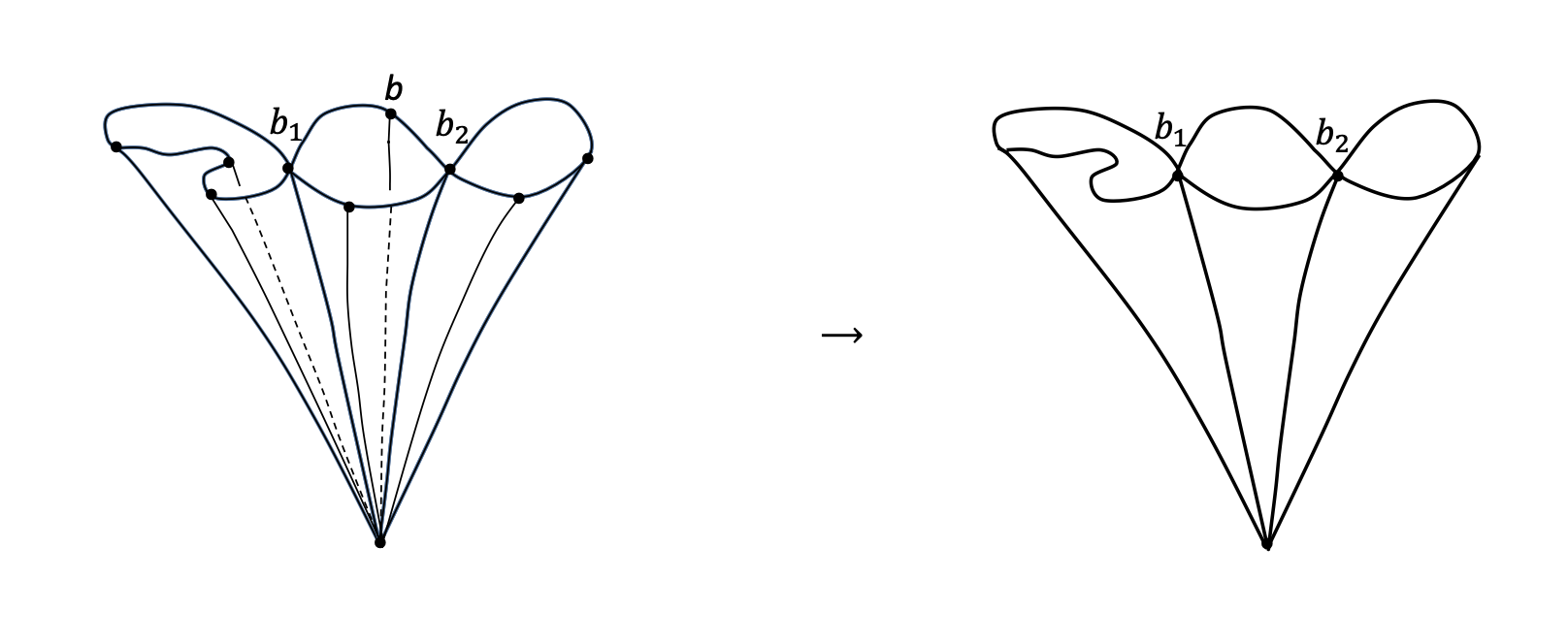}}
\caption{Simplification process.}
\end{figure}

\begin{theorem}[Birbrair Classification Theorem, Theorem 8.1 \cite{Lev}] \label{Clas_theo}
{\rm 
Let $X_1, X_2 \subset \R^p$ be two dimensional semialgebraic subsets with $a_1 \in X_1$  and $a_2 \in X_2$. The germs $(X_1, a_1)$ and $(X_2, a_2)$ are (abstract) bi-Lipschitz inner  equivalent if and only if the Canonical H\"older Complexes of $X_1$ at $a_1$ and $X_2$ at $a_2$ are combinatorially equivalent. 

}
\end{theorem}


\subsection{Bi-Lipschitz classification in the case $(n,p) = (2,3)$.}

Let $f,g : (\R^2,0) \to (\R^3,0)$ be map germs. 
We denote by $X_f = f(U)$ and $X_g = g(U)$, where $U$ is a sufficiently small neighbourhood of $0$. 

Let $\pi: \R^3 \to \R^2$ be a linear projection which is generic for $X_f$. That means a finite map such that 
the topological closure of $\pi\vert_{(X_f)_{reg}}$  is the polar curve of $f$ and we denote it by $\Gamma^1(X_f)$.

As an application of  Birbrair's Theorem, we have the following proposition: 
\begin{proposition} \label{demai}
Let $f,g : (\R^2,0) \to (\R^3,0)$ be finitely determined map germs. Then  $X_f$ and $X_g$ are (abstract) bi-Lipschitz inner equivalent  if and only if 
the following invariants are bi-Lipschitz equivalent:
\begin{enumerate}
\item $link(f)$ and $link (g)$; 
\item $f(D(f))$ and $g(D(g))$;
\item the polar curves $\Gamma^1(X_f)$ and $\Gamma^1(X_g)$.
\end{enumerate}

\end{proposition}

\begin{preuve}
As $f, g:(\R^2,0) \to (\R^3,0)$ are finitely determined map germs, one can consider $f$ and $g$ as polynomial maps. Then $X_f $ and $X_g$ are  closed semialgebraic surfaces. 
By Theorem \ref{theorem4.7}, there exists a semialgebraic H\"older Complex $(\Gamma_f, \beta)$ such that $X_f$ is a semialgebraic Geometric H\"older Complex corresponding to $(\Gamma_f, \beta)$ 
with the principal vertex $0$. Similar property holds for $X_g.$
Using Birbrair's simplification process, the germs 
$X_f$ and $X_g$ correspond to 
two Canonical H\"older Complexes at 0, that we denote by 
$\tilde{\Gamma}_f$ and $\tilde{\Gamma}_g$. Let us suppose that the invariants in (1), (2) and (3) are bi-Lipschitz equivalent. Then, since $link(f)$ and $link(g)$ are bi-Lipschitz equivalent, the graphs $\Gamma_f$ and $\Gamma_g$ are homeomorphic. 
 Then the loop vertices and the critical vertices of $\Gamma_f$ correspond res\-pec\-tive\-ly to the loop vertices and  the critical vertices of $\Gamma_g$. 
Notice that the critical vertices of $\Gamma_f$ 
are precisely the intersection points of $link(f)$ and $f(D(f)).$ 
 We have the same thing for $\Gamma_g$. 
 Then the branches of $f(D(f))$ correspond to the branches of $g(D(g))$. 
 By Theorem \ref{LA}, since $f(D(f))$ and $g(D(g))$ are bi-Lipschitz equivalent then the contact order $\beta$ of two branches of $f(D(f))$ is the same the contact order of the two corresponding branches of $g(D(g))$. 
 Now, by hypothesis, we can choose the generic projection $\pi$ in such a way that each branch of the polar curve 
 corresponds to a loop vertex, for $f$ and $g$.
Hence, the two Canonical H\"older Complexes $\tilde{\Gamma}_f$ and $\tilde{\Gamma}_g$ are combinatorially equivalent. 
The proof of the Proposition follows from Birbrair's Classification Theorem \ref{Clas_theo}. 

For the converse, let us suppose that $X_f$ and $X_g$ are semialgebraically bi-Lipschitz equivalent with respect to the inner metric. From the Classification Theorem \ref{Clas_theo}, it follows that the associated H\"older Complexes $\tilde{\Gamma}_f$ and $\tilde{\Gamma}_g$ are combinatorially equivalent. 
As a consequence, we get the outer bi-Lipschitz equivalence between $f(D(f))$ and $g(D(g))$ 
and also bi-Lipschitz equivalence between the polar curves 
$\Gamma^1 (X_f)$ and $\Gamma^1 (X_g)$. Notice that this bi-Lipschitz equivalence can be extended to the ambient space (see  \cite{LA}).  Furthermore,  as $f$ and $g$ are finitely determined, the only singularities of $link(f)$ and $link(g)$ are transverse double points, then  these two sets are Lipschitz normally embedded. Hence, it follows that $link(f)$ and $link(g)$ are outer bi-Lipschitz equivalent. 

\end{preuve}

\begin{example}  \label{counter-example}
In this example, we provide real analytic map-germs for which $D(f) = \{ 0\}$. 
 Let us consider $f, g: (\R^2, 0) \to (\R^3, 0)$ defined by 
$$f(x, y) = (x, x^2y+y^3, y^5), \qquad g(x,y) = (x, x^4y+y^5, y^7).$$ 
Let $X_f$ and $X_g$ be the images $f(\R^2)$ and $g(\R^2),$ respectively.
In this case $\Sigma(f) = \Sigma(g) = \{(0,0)\}$. 
For determining the double point set $\widetilde{D}(f)$, we solve the following system 
$$
\begin{cases}
\dfrac{x^2y + y^3 - x^2y' - y'^3}{y-y'} = 0 \cr 
\dfrac{y^5 - y'^5}{y - y'} = 0.
\end{cases}
$$
The first equation implies $x^2 + p^2(y, y') = 0$ 
and the second equation can be written as $ q^4(y, y') = 0$, where $p^2$ and $q^4$ are the 
symmetric polynomials of degree 2 and 4, respectively. 
We can prove that there exist constant numbers $C$ and $\tilde{C}$ such that 
$$C\, (y^4+y'^4) \leq  q^4(y, y') \leq \tilde{C}\, (y^4+y'^4)$$
in a small neighbourhood of the origin. 
This shows that the only solution of the second equation is $y = y' =0$. 
Consequently, 
$$\widetilde{D}(f) = \{x = y = y' = 0\}.$$
With the same arguments, we can also show that 
$$\widetilde{D}(g) = \{x = y = y' = 0\}.$$

\vspace{0.1cm}
We can see that $f$ and $g$ satisfy the following conditions:

(1) $f$ and $g$ are finitely determined,

(2) The  double point sets $D(f)$ and $D(g)$ as well as their images 
$f(D(f))$ and $g(D(g))$ are reduced to the origin, {\it i.e.} $f(D(f))=0$ and $g(D(g))=0$,

(3) The generic polar curves $\overline{\Sigma \pi_f \vert_{(X_f)_{reg}}}$ and 
$\overline{\Sigma \pi_g \vert_{(X_g)_{reg}}}$ where $\pi_f : \R^3 \to \R^2 $ and 
$\pi_g: \R^3 \to \R^2 $ are generic linear projections, are also reduced to the origin,


(4) the hyperplanes sections of $X_f$ and $X_g$ are all smooth, excepted for $x=0$, which are respectively $(y^3, y^5)$ and $(y^5, y^7)$, 

(5) The equations of ${X_f}$ and ${X_g}$ are respectively 
$$z^3 -y^5 + 5 x^2yz^2 + 5 x^4y^2z + x^{10}z = 0$$
$$ z^5 - y^7 + 7 x^4y^2z^3 + 7 x^8y^4z - 7 x^{16}y z^2 + x^{28}z = 0.$$

Then  ${X_f}$ and ${X_g}$ are not ``outer'' biLipschitz equivalent due to the fact that 
their multiplicities are different. 

\vspace{0.2cm}

However, as the tangent cone of $X_f$ and $X_g$ are 2 dimensional planes, it follows from  Birbrair Classification Theorem that
$X_f$ and $X_g$ are bi-Lipschitz equivalent to a plane. Then, $X_f$ and $X_g$ are inner bi-Lipschitz equivalent.

\vspace{0.2cm} Notice that we can  generalize this example in order to obtain an infinite family
$$f_k(x,y)= (x, x^{2k}y+y^{2k+1}, y^{2k+3}), \, k \geq 1$$ of $\mathcal A$-finitely determined map-germs, such that the
images $f_k(\mathbb R^2)=X_{f_k}$ are all inner bi-Lipschitz equivalent, however $f_k$ and $f_k'$ are not outer bi-Lipschitz equivalent for $k \neq k'.$

\begin{figure}[H]
\scalebox{0.1}{\includegraphics{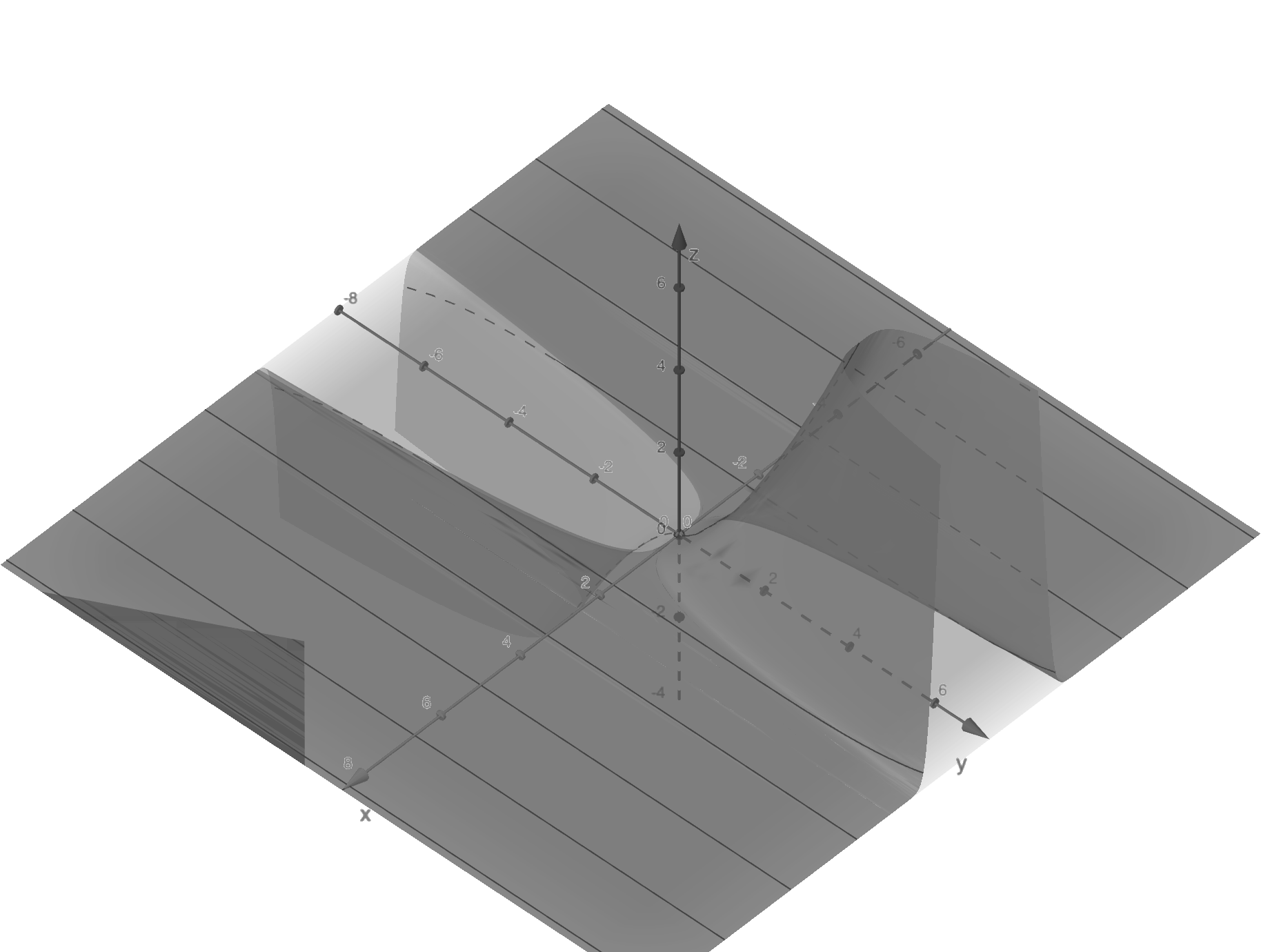}}
\qquad \qquad 
\scalebox{0.1}{\includegraphics{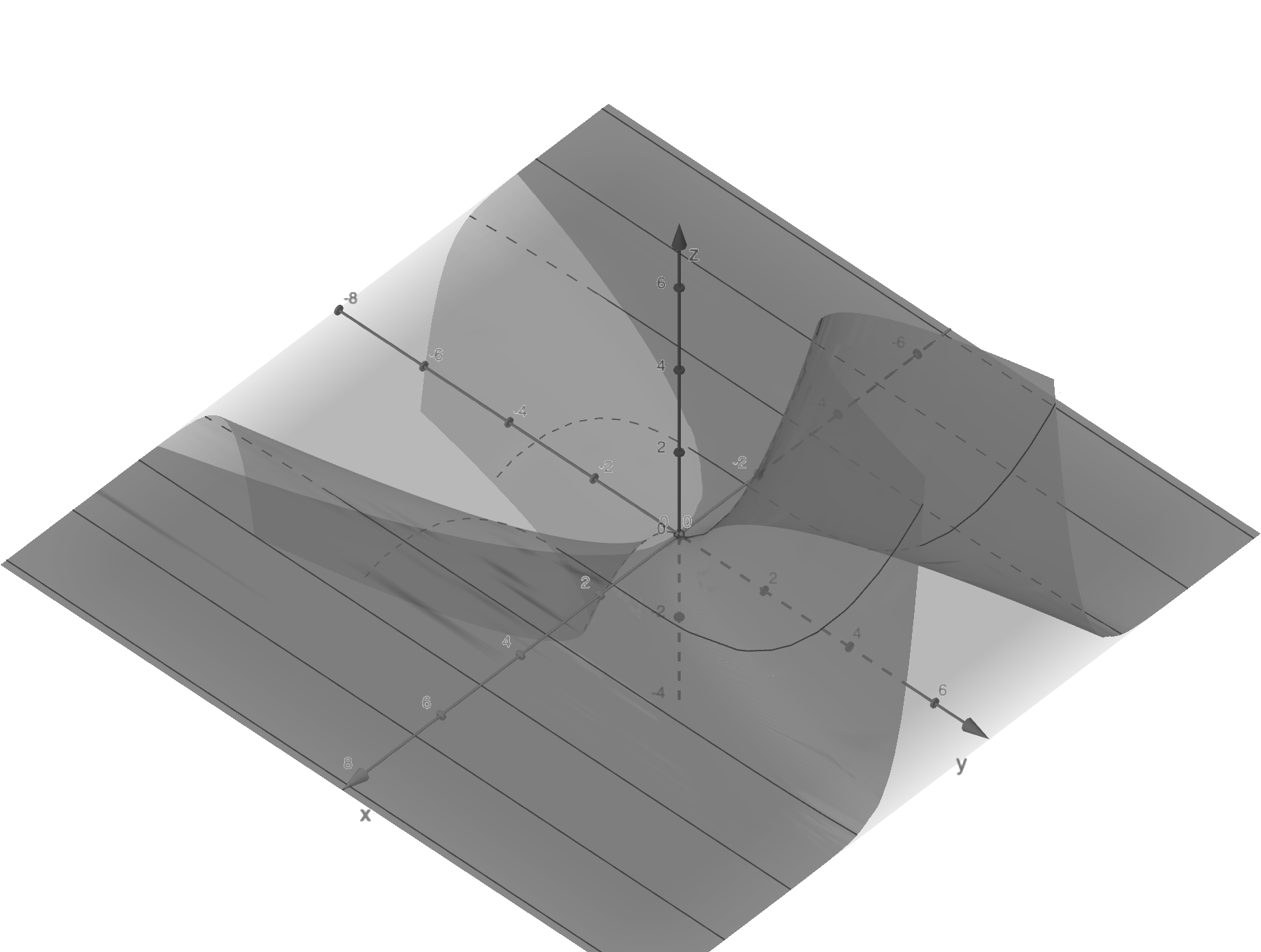}}
\caption{$X_f$ and $X_g$.}
\end{figure}
\end{example}

\begin{corollary}
Let $f,g : (\R^2,0) \to (\R^3,0)$ be finitely determined map germs such that the following invariants are  bi-Lipschitz equivalent:
\begin{enumerate}
\item $link(f)$ and $link (g)$; 
\item $f(D(f))$ and $g(D(g))$;
\item the polar curves $\Gamma^1(X_f)$ and $\Gamma^1(X_g)$. 
\end{enumerate}
Then if $X_f$ and $X_g$ are LNE, the surfaces $X_f$ and $X_g$ are bi-Lipschitz equivalent 
with respect to the outer metric.
\end{corollary}

\section{Surfaces with homogeneous parametrization} \label{BiLipschitzforquasihomogenousSurfaces}

Let $f, g: (\R^2, 0) \to (\R^3, 0)$ be finitely determined map germs of corank 1 with  homogeneous 
parametrization of the same degrees,  that means that $f$ and $g$ can be written in the following way
\begin{equationth} \label{form_homogeneous_parametrization}
f(x, y) = (x, f_{2}(x,y), f_{3}(x,y)), \quad g(x, y) = (x, g_{2}(x,y), g_{3}(x,y))
\end{equationth}
where $f_j$ and $g_j$ are homogeneous polynomials of the same degree $d_j$, for $j =2, 3,$ 
and gcd$(d_2,d_3)=1 \,  \text{or}\, 2$ (See \cite{Ot} , Lemma 6.2).

Notice that the equation $\phi(x,y,z) =0$ defining $X_f$ ({\it resp.}, $X_g$) is not homogeneous in general. 

The following theorem shows that in this case, in order to obtain the bi-Lipschitz e\-qui\-va\-lence 
of $X_f$ and $X_g$, 
we do not need any hypothesis on the bi-Lipschitz equivalence of the set $f(D(f))$ and $g(D(g))$.

The proof of the theorem uses the notion of {\sl identification component } 
due to Ruas and Silva (\cite{Ot},  Definition 2.4): 
Take a representative $f:U\to V$ of $f$, where $U$ and $V$  are neighbourhoods of 
$0$ in $\R^2$ and $\R^3$.  An {\sl identification component} of $D(f)$ is 
an irreducible component $D(f)^j$ of $D(f)$ such that 
the restriction $f\vert_{D(f)^j}: D(f)^j \to \R^3$  is generically one to one. 
If the restriction is generically $2 - 1$, we say that $D(f)^j$is a
{\sl fold component} of $D(f)$.

\begin{theorem} \label{theo_quasehomogeneo}
Let $f, g: (\R^2, 0) \to (\R^3, 0)$ be finitely determined map germs of corank 1 with homogeneous parametrization of the form (\ref{form_homogeneous_parametrization}).  
 If $link(f)$ and $link(g)$ are bi-Lipschitz equivalent  and the polar curves $\Gamma^1(X_f)$ 
  and $\Gamma^1(X_g)$ are bi-Lipschitz equivalent, 
 then  $X_f$ and $X_g$ are bi-Lipschitz equivalent with respect to the inner metric.
\end{theorem}

\begin{preuve}
Let us write $f$ in the form (\ref{form_homogeneous_parametrization}), that means 
$$f(x, y) = (x, f_2(x,y), f_3(x,y)).$$
According to \cite{MM}, the set of double points $\tilde D(f)$ is the set of the solutions of the following system:
$$\frac{f_2(x,y) - f_2(x,y')}{y-y'} = 0, \quad \frac{f_3(x,y) - f_3(x,y')}{y-y'} = 0.$$ 

\noindent Then, since $f_2$ and $f_3$ are homogeneous polynomials,  $D(f) \subset \mathbb R^2$ is the union of lines $\gamma_i$ whose equations are of the type $y = \alpha_i x$ and in some cases also the line $x =0$. In  the case $y = \alpha_i x$,  we have 
$$f(\gamma_i) = (x, f_2(x, \alpha_i x), f_3(x, \alpha_ix)).$$
We see that $f_2(x, \alpha_i x)$ is  a homogeneous polynomial of degree $d_2 $  and $f_3(x, \alpha_i x)$ 
is  a homogeneous polynomial of degree $d_3$ with respect to the variable $x$.  When  gcd{$(d_2, d_3)=1$}, it follows
that $f_| {\gamma_i}$ is injective and each branch of $f(D(f)),$ (respectively of $g(D(g))$) is an {\it identification component} of two branches
say, $\gamma_i$ and $\gamma_j$  (see  \cite{Ot}, Section 6 for details).
In this case,  the image of the set of double points $f(D(f))$ is the union of 
parametrized curves of the form $(x, \tilde{f}_2, \tilde{f}_3)$ where $\tilde{f}_j$ are homogeneous 
polynomials of degree $d_j$ with respect to variable $x$, for $j =2, 3$. 
Similarly, $g(D(g))$ is the union of parametrized curves  of the form $(x, \tilde{g}_2, \tilde{g}_3)$ where 
$\tilde{g}_j$ are homogeneous polynomials of the same degree than $\tilde{f}_j$, for $j =2, 3$. 
Therefore, according to Theorem \ref{LA}, 
the order of contact of two branches of $f(D(f))$ is the same than the order of contact of 
two corresponding branches of  $g(D(g))$. 

In the case that gcd{$(d_2,d_3)=2,$}  the equation  $x = 0$ defines a branch of $D(f)$ and $D(g)$, it provides a {\it fold} component,   respectively of $f(D(f))$ and $g(D(g))$ (see Proposition 6.2, \cite{Ot}).

Now, with the hypothesis that $link(f)$ and $link(g)$ are bi-Lipschitz equivalent and the polar curves $\Gamma^1(X_f)$ 
  and $\Gamma^1(X_g)$ are bi-Lipschitz equivalent, the proof follows from Proposition \ref{demai}.

\end{preuve}

\begin{corollary} \label{corquasehomo}
Let $f, g: (\R^2, 0) \to (\R^3, 0)$ be finitely determined map germs of corank 1 
with homogeneous parametrization such that $link(f)$ and $link(g)$ are bi-Lipschitz equivalent and the polar curves $\Gamma^1(X_f)$ 
  and $\Gamma^1(X_g)$ are bi-Lipschitz equivalent. 
  Then if $X_f$ and $X_g$ are LNE,   
the surfaces  $X_f$ and $X_g$ are bi-Lipschitz equivalent with respect to the outer  metric.
\end{corollary}

In the next proposition we  apply the previous results to relate the $C^{0}- \mathcal A$ classes of finitely determined map germs $f$ of corank 1 with  homogeneous  parametrization 
 and the inner bi-Lipschitz type of $X_f.$

First we introduce some notation. Let $H_{n,m}$ be the  space of real polynomial mappings
$$ f(x,y):= (x, p(x,y), q(x,y)),$$ where $p(x,y)=\sum^{n}_{i=0}a_ix^{n-i}y^i \,\,\, and \,\, q(x,y)=\sum^{m}_{i=0}b_ix^{m-i}y^i,$  with $\text{gcd}(n,m)=1, \,\, n< m.$

We can write 
$$ p(x,y)=y^n+x\tilde p(x,y), \,\,\,\,\, q(x,y)=b y^m+x\tilde q(x,y),$$ where $\tilde p(x,y), \,\, \tilde q(x,y)$ are homogeneous polynomials of degrees $n-1$ and $m-1,$  respectively.

Let $\pi_{\beta} : \mathbb R^3 \to \mathbb R^2$ be the linear projection $(X,Y,Z) \mapsto (X, Y+\beta Z), \, \beta \in \mathbb R.$

We say that $\pi_{\beta}$ is generic for $f \in H_{n,m}$ if the map-germ 
$$ h= \pi_{\beta} \circ f: (\mathbb R^2,0) \to (\mathbb R^2,0)$$  is finite and its singular set
$\Sigma(h)$ is a reduced plane curve $C(f).$ Notice that  $C(f)$ is defined by the equation $p_y+\beta q_y=0.$
The polar curve $\Gamma^1(f)$ is the space curve $f(C(f)).$

\begin{proposition} \label{Proposition 5.4}
	With the above notation, the following conditions hold:
	\begin{itemize}
\item [(1)]	The subset $F_{n,m}:=\{f\in H_{n,m}\,| \,\, f \,\,\text{is}\,\, \mathcal A$-$\text{finitely determined} \, \}$ is a non empty open and dense subset of $H_{n,m}.$

\item [(2)] The $C^{0}-\mathcal A$-type of $f \in F_{n,m}$ is determined by the Gauss word of $link(f).$

\item [(3)] The Gauss word of $link(f)$ determines the critical vertices and the loop vertices of
the canonical H\"older Complex of $X_f.$ Moreover, if $c_v$ and $l_v$ denote the number of critical vertices and loop vertices, respectively, we have that $$c_v \leq \, \frac{nm-(n+m)+1}{2} \, \, \text{and} \,\, c_v \cong \frac{nm-(n+m)+1}{2}\, (\text{mod}(2)),$$ and $l_v$ is equal to the number of syllables of type $a_j a_j^{-1}$ and $a_j^{-1}a_j$ in the Gauss word of $link(f).$

\item [(4)] If $f, g \in F_{n,m},$ then 
$X_f$ and $X_g$ are inner bi-Lipschitz  equivalent if and only if $f$ and $g$  are ${\mathcal{C}^0}$-$\mathcal{A}$-equivalent and the curves $C(f)$ and $C(g)$ are outer 
bi-Lipschitz equivalent.

\end{itemize}

\end{proposition}

\begin{preuve}

(1)	and (2): We first notice that $F_{n,m}\neq \emptyset$ since $f(x,y)=(x, y^n, (x+y)^m),$ where  $gcd(n,m)=1$, is finitely determined (see \cite{Gui}, Proposition 9.8).
Then $F_{n,m}$ is an open and dense subset of $H_{n,m}$ follows because finite
determinacy is an open condition (\cite{Damon}) and its complement is a proper semialgebraic set in $H_{n,m}$
	{\rm(see \cite{Gui} and \cite{Ot}, Example 5.5 and Lemma 7.1)}.
	 Moreover, on each connected component of $F_{n,m}$, 
	 the $C^{0}$-$\mathcal{A}$ type is constant 
	  (see Damon's Theorem 1 in \cite{Damon}, page 381). 
	 It follows from Corollary \ref{theorem-thon-juan}  that the $C^0$-$\mathcal A$ type of $f$ is determined by
	  the Gauss word of $link(f).$
	  
	  To prove (3), notice that the critical vertices of the canonical H\"older Complex of $X_f$ 
	  are given by the double points of the link of $X_f$ 
	  and that each syllable $a_ja_j^{-1}$ or $a_j^{-1}a_j$ determines a loop in this complex.
	  
	  Now, by Proposition 6.2 in  \cite{Ot}  it follows that the double point set $D(f_{\mathbb C} )$  of the complexification of $f$  is the germ of a homogeneous curve with $d$ smooth irreducible components, where $d=nm-n-m+1$ and 
	  
	  \begin{center}
	  	$D(f_{\mathbb C})= V( \displaystyle {  \prod_{i=1}^{d}}(x-\alpha_iy))$
	  \end{center}
	  
	  \noindent where $\alpha_i \in \mathbb{C}$.
	  
	  Then, the set $D(f) \subset \mathbb R^2$ is a union of lines through the origin, and the number of lines
	  is less than or equal to $d, $ and is congruent to $d$ module $2.$ Since $gcd(n,m)=1,$ then all lines are identification lines
	  so that $c_v\cong  \, \frac{d}{2} \, (mod 2).$  The condition on the loop vertices follow trivially.

	 (4) It follows from  Theorem \ref{theo_quasehomogeneo} that for $f, g \in F_{n,m},$
	 $X_f $ and $X_g$ are inner bi-Lipschitz equivalent if and only if $link(f)$ is bi-Lipschitz equivalent to
	 $link(g)$ and $\Gamma^1(f) $ is bi-Lipschitz equivalent to $\Gamma^1(g).$ The result now follows
	 from  Proposition \ref{link-cond} and (3), noticing that 
	 $C(f)=f^{-1}(\Gamma^1(f))$ and $C(g)=g^{-1}(\Gamma^1 (g)).$

\end{preuve}

\vspace{0.2cm}

\noindent {\bf Question:} Given $f, g \in F_{n,m},$ does it follow that 
$$f \backsimeq_{C^0-\mathcal A} g  \Longleftrightarrow X_f \, \text {and } X_g$$ are inner bi-Lipschitz 
equivalent?

\vspace{0.2cm}

\noindent{\bf Acknowledgement:} We thank Otoniel Nogueira da Silva, Alexandre Fernandes 
and Edson Sampaio 
for many useful conversations. The first named author thanks, on the Brazilian side, the Project UNESP-FAPESP n${}^{\rm o}$ 
2015/06697-9, the Project USP-COFECUB n${}^{\rm o}$  UC Ma 163-17 and the Project USP-FAPESP n${}^{\rm o}$ 2018/07040-1, and, on the French side, CNRS (UMR 7373) and Aix-Marseille University (I2M).
The second named author was partially supported by FAPESP, grant $\# 2019/21181-0$ and CNPq, grant $\# 433257/2016-4$ and by CNPq, grant $\#  305695/2019-3.$ 
The third named author was supported by the Brazilian S\~ao Paulo Research Foundation (FAPESP) 
grant  \# 2023/07802-7. 
\bibliographystyle{plain}

\begin{thebibliography}{BBFGK}

\def\rond{\mathaccent"7017}
\def\omini{\raise 1ex \hbox{\sevenrm o}}


\bibitem{Lev} L. Birbrair, {\it Local bi-Lipschitz classification of 2-dimensional semialgebraic sets}, 
Houston Journal of Mathematics (1998), vol. 25, no. 3, 453--472.


\bibitem{LA} L. Birbrair and A. Fernandes, {\it Metric theory of semialgebraic curves}, Revista Matem\'atica Complutense (2000), vol. XIII, num. 2, 369--382.


\bibitem{Damon} J. Damon, {\it Topological triviality and versality for subgroups of A and K: II. Sufficient conditions and applications}, Nonlinearity, v.5, n.2, p. 373-412. 1992.

\bibitem{DPW} A.A. Du Plessis and C.T.C. Wall, {\it Generic projections in semi-nice dimensions}, Compositio Mathematica, (2003), Vol 135, 179--209.


\bibitem{Fukuda} T. Fukuda, {\it Local topological properties of differentiable mappings I}, Invent. Math. 65 (1981/1982) 227-250.

\bibitem{Gaf} T. Gaffney, {\it Polar multiplicities and equisingularity of map germs}, Topology vol.32 No.1 (1993) 185--223.

\bibitem{MM} W.L. Marar and D. Mond, {\it Multiple point schemes for corank 1 maps}, J. London Math. Soc., 39 (1989), 553--567.

\bibitem{Thon1} W.L. Marar and J.J. Nu\~no-Ballesteros, {\it The doodle of a finitely determined map germs from $\R^2$ to $\R^3$}, Advances in Mathematics {\bf 221} (2009), pages 1281-1301. 

\bibitem{Thon2} W.L. Marar, J.J. Nu\~no-Ballesteros and G. Pe\~nafort-Sanchis, {\it Double point curves for corank 2 map germs from $\C^2$ to $\C^3$}, Topology and its Applications {\bf 159} (2012), pages 526-536. 

\bibitem{Mat} J. Mather, {\it Stability of $C^\infty$ -mappings. VI. The nice dimensions}, Proc. Liverpool Singularities-Sympos., I (1969/70), Lecture notes in Math., vol. 192, Springer-Verlag, Berlin and New York, 1971, pp. 207-253. 

\bibitem{MaSa} R. Mendes, J.E. Sampaio, {\it On the link of Lipschitz normally embedded sets}, International Mathematic Research Notes, pp 1-14, 2023.

\bibitem{Mond1} D. Mond, {\it Some remarks on the geometry and classification of germs of map from surfaces to 3-space}, Topology 26, (1987), 3, 361-383.

\bibitem{MNB} D. Mond and J. Nu\~no-Ballesteros, {\it Singularities of mappings} 2016, available on  \hfill\\homepages.warwick.ac.uk/~masbm/LectureNotes/book.pdf

\bibitem{NP} W. Neumann, A. Pichon, {\it Lipschitz geometry of complex surfaces : analytic invariants and equisingularity.} arXiv:1211.4897, 2012.

\bibitem{NRT} N. Nguyen, M. Ruas, S. Trivedi, {\it Classification of Lipschitz simple function germs.} 
Proc. London Math. Soc., 121(1), 51--82, 2020. 

\bibitem{NR} J. Nu\~no-Ballesteros and R. Mendes, {\it Topological classification and finite determinacy of knotted maps},  ArXiv: 1811.01113.

\bibitem{Gui} G. Pe\~{n}afort Sanchis, {\it Reflection maps}, Mathematische Annalen, 2020. 

\bibitem{PP1} A. Parusi\'nski, L. Paunescu,  {\it Lipschitz Stratification of Complex Hypersurfaces in Codimension 2,} J. Eur. Math. Soc. 25, 1743--1781 (2023).

\bibitem{PP2} A. Parusi\'nski, L. Paunescu,  {\it Zariski's dimensionality type of singularities. Case of dimensionality type 2.} arXiv:2104.07156 [math.AG].

\bibitem{Cidinha} M.A.S. Ruas, {\it Basics on Lipschitz Geometry},  p. 111-155, in: {\bf Introduction to Lipschitz Geometry of Singularities}, Lecture Notes of the International School on Singularity Theory and Lipschitz Geometry, Cuernavaca, June 2018.

\bibitem{Ot} M.A.S. Ruas and O.N. Silva, 
{\it Whitney equisingularity of families of surfaces in $\C^3$},  
Math. Proc. Cambridge Philos. Soc. 166 (2019), no. 2, 353-369.



\bibitem{Ot1} O.N. Silva, {\it Surfaces with non-isolated singularities}, Thesis USP, S\~ao Carlos, 2017, avalaible on 
 http://www.teses.usp.br/teses/disponiveis/55/55135/tde-10052017-085440/pt-br.php

\bibitem{Wall} C.T.C. Wall, {\it Finite determinancy of smooth map-germs}, Bull. London Math. Soc., {\bf 13} (1981), 481-539. 

\bibitem{Whitney} H. Whitney, {\it The singularities of a smooth $n$-manifold in $(2n-1)$-space}, Ann. of Math. {\bf 2}, 45 (1944), 247-293. 

\end{thebibliography}

\end{document}